# Efficient solutions for nonlocal diffusion problems via boundary-adapted spectral methods

Siavash Jafarzadeh[1], Adam Larios[2,*], and Florin Bobaru[1,*]

[1] *Department of Mechanical and Materials Engineering, University of Nebraska-Lincoln, Lincoln, NE 68588-0526, USA*

[2] *Department of Mathematics, University of Nebraska-Lincoln, Lincoln, NE 68588-0130, USA*

*corresponding authors. Email addresses: alarios@unl.edu (A. Larios); fbobaru2@unl.edu (F. Bobaru)

**Abstract**

We introduce an efficient boundary-adapted spectral method for peridynamic diffusion problems with arbitrary boundary conditions. The spectral approach transforms the convolution integral in the peridynamic formulation into a multiplication in the Fourier space, resulting in computations that scale as $O(N\log N)$. The limitation of regular spectral methods to periodic problems is eliminated using the volume penalization method. We show that arbitrary boundary conditions or volume constraints can be enforced in this way to achieve high levels of accuracy. To test the performance of our approach we compare the computational results with analytical solutions of the nonlocal problem. The performance is tested with convergence studies in terms of nodal discretization and the size of the penalization parameter in problems with Dirichlet and Neumann boundary conditions.

1. Introduction

Nonlocal models have been introduced to address certain phenomena which local models fail to describe satisfactorily. Delayed reaction-diffusion in biology [1], swarm of organisms [2], pedestrian traffic [3], flocking of birds [4-6], plane waves in solids [7], elasticity of nano-beams [8], and material damage [9,10] are some examples of problems where nonlocal models are useful. Material damage models in particular are of significant interest, being used for failure prediction of critical materials and structures. Physical features in damage (evolving cracks and distributed failure) and small-scale heterogeneities can be naturally modeled using nonlocal approaches [11,12], and would be otherwise difficult to describe or prohibitively expensive to compute with classical local approaches. Peridynamics, as a nonlocal extension of continuum mechanics [13,14], has been successful in modeling damage evolution and material failure [15,16,13]. Dynamic brittle fracture [17-19], fatigue and thermally-induced cracking [20,16], fracture in porous and granular materials [21-23], failure of composites [24,25], corrosion damage [26-29], and stress corrosion cracking [30-32], are among some applications of this formulation in modeling material damage.

In peridynamics (PD), material behavior at each point $x$ depends on the interactions of that point with all of the points $\hat{x}$ in its neighborhood [14]. This neighborhood (usually a line segment in 1D, a disk in 2D, and a sphere in 3D) centered at $x$ is called the "horizon region" of $x$ and is denoted by $H_x$. $H_x$ is the subdomain where nonlocal interactions exist for $x$.

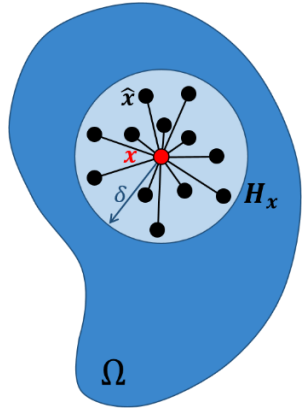

Fig. 1. Nonlocal interactions of point $x$ with its neighboring points $\hat{x}$, in its horizon region $H_x$, in a schematic peridynamic body $\Omega$.

Mathematical models of physical behavior using this approach are in the form of integro-differential equations, as spatial derivatives in the classical PDEs are replaced by convolution integrals that integrate the pairwise interactions of $x$ with the points in $H_x$. Integration has significantly more relaxed smoothness and continuity requirements compared with differentiation and hence, it allows for more robust handling of discontinuities, such as cracks. While nonlocality facilitates describing material degradation and provides certain advantages for incorporating small-scale features into large-scale models [33,34], it also adds a significant computational cost, due to the convolution integral involved, compared with local models.

Two types of numerical methods have been commonly used for the discretization of PD models. One popular method that offers much flexibility for arbitrary/unguided damage/fracture evolution is a meshfree method based on the one-point Gaussian quadrature of the integral operator [35,36]. If the total number of nodes in the domain is $N$ and the number of nodes inside the horizon of each point is $M$, the computational cost at each time step in an explicit algorithm will be $O(NM)$. Note that in one-dimension, $= \frac{\delta}{L} N$, where $L$ is the length of the domain and $\delta$ is the horizon size (see Fig. 1). Therefore, for a fixed horizon size, $M$ itself varies as $O(N)$, and the computational cost is $O(N^2)$.

Finite element (FE) methods have also been used to discretize PD models: in some of such models, each pairwise interaction (bond) is represented as a truss element [37,38], while others use continuous or Discontinuous Galerkin (DG) discretization methods, for example [39-41]. In all FEM-based discretizations of PD models, explicit solutions also cost $O(N^2)$ per time step. We note that regular FEM discretizations are not used for modeling of problems in which discontinuities appear due to the inherent difficulties of the method (see [39]). This is the main reason that only truss-based or DG methods have been used in PD models of failure/fracture, in addition to the most successful, meshfree discretization.

Coupling local models (discretized, e.g., with FEM) with PD models (discretized, e.g., with the meshfree method) has been seen as one way to increase the efficiency of simulations with PD models [42,43]. These approaches are only beneficial when the region where nonlocality is helpful or dominates, covers only a small portion of the system modeled (e.g. small localized damage or crack growth). The advantage is lost in problems in which, for example, failure is affecting a large part of the domain [15,19].

Note that in the methods mentioned above for discretizing nonlocal models, $M$ increase exponentially with the problem dimension. Indeed, assume that the length scale of the domain is $L$ and the grid spacing in each direction is $\Delta x$. If $\delta$ is the radius of the neighborhood (also called the "horizon size", or simply "the horizon"), then $N = \left(\frac{L}{\Delta x}\right)^d$ and $M \propto \left(\frac{\delta}{\Delta x}\right)^d$ where $d$ is the dimension of the problem. Computational cost per times step can then be expressed as $\left[\left(\frac{L}{\Delta x}\right)^d \cdot \left(\frac{\delta}{\Delta x}\right)^d\right]$.

In Fourier spectral methods, the solution is transformed to Fourier space (if the solution is assumed *periodic*), and the governing equations is reformulated based on the transformed solution. In the case of classical PDEs, spatial derivatives transform to multiplication operators, and the PDE reduces to a system of ODEs in Fourier space, which is far easier to solve [44]. For nonlocal models, Fourier transformation disentangles the convolution integral and reduces it to a multiplication in the spectral space. The only major cost in the spectral method is the Fourier transform itself, and its inverse. For this, the well-known *Fast Fourier Transform* (FFT) algorithms are available, at a cost of $O(N\log N)$ [45,46,44]. Not only is the cost of the Fourier spectral method significantly lower than the two other numerical methods used to discretize PD models, but the FFT is also easily parallelized, further increasing the potential advantages of this approach.

Although the spectral method seems to be a promising candidate for computing solutions to nonlocal problems, the assumption of periodicity limits its application. Most real-world problems are not periodic. A few recent studies have introduced Fourier spectral methods for periodic nonlocal models. For example, this method has been used for the nonlocal Allen-Cahn equation [47], nonlocal damage models [48], and peridynamic nonlocal operators for diffusion and wave propagation problems [49-51]. In all these cases, the problems considered were periodic. Another method, while not spectral, uses the FFT to diagonalize the stiffness matrixes arising in FE and collocation discretization methods for nonlocal problems with non-periodic boundary conditions [52,53]. Although the order of computation is $N\log N$, the method is restricted to simple domain shapes like a square in 2D. The method is also dependent to the horizon shape. The authors of [52] state that the method is not applicable to domains with complex geometry, or to heterogeneous domains, and is challenging to use in 3D.

While the periodicity of the solution is inherent in Fourier spectral methods, there exist ways to overcome this limitation and apply them to general problems with non-periodic boundary conditions. *Volume penalization* is one of such method.

Penalization methods have been used with the local Navier-Stokes PDE to introduce solid obstacles/boundaries in fluid flows, without changing the equations and discretization. A rigorous, simple volume-penalization method based on Brinkman model for flow in porous media [54] is developed by Angot et al. [55]. In [55], a large viscous term is added to the equation in the solid region to impose a Dirichlet (no-slip) boundary condition for the fluid-solid contact. Kevlahan and Ghidaglia [56], used this method with the Fourier spectral method in fluid dynamics problems. In these methods, the solution (velocity) is penalized by a substantially higher viscosity in the solid region to enforce a zero-velocity boundary condition. The method was applied for modeling flow over stationary or moving solid obstacles with complex geometries, inside periodic or confined fluid domains [57-59]. Volume penalization has also been used to enforce no-flux (Neumann-type) boundary conditions in advection-diffusion problems solved with the spectral method [60]. Another example for this efficient method is the 3D simulation of bumblebee flight in wind flow [61,62].

In the present study, we introduce a spectral method to obtain efficient solutions to nonlocal equations of the peridynamic type for diffusion with arbitrary, non-periodic boundary conditions, using the volume-penalization technique. In Section 2, the PD formulation and boundary conditions implementation in PD problems are briefly discussed. Spectral methods and volume penalization for PD problems are introduced in Section 3. Stability analysis is provided in Section 4 and two examples with non-periodic BCs are solved in Section 5. Convergence studies are provided in Section 6.

## 2. Peridynamic nonlocal formulation

We start our development of spectral methods for peridynamic models with the PD diffusion equation in 1D. The methods described here are, however, applicable to other PD models, as well as to any other model with convolution integrals. Eq. (1) is the general form of the PD diffusion equation in 1D [63]:

$$\frac{\partial u(x,t)}{\partial t} = \nu \mathcal{L}_\delta u(x,t) + f(x,t) \tag{1}$$

where $x$ is the position in the 1D domain $\Omega$, $u(x,t)$ is the unknown (the solution field) at point $x$ and time $t$, $\nu$ is the diffusivity, $\mathcal{L}_\delta$ is the PD Laplacian operator (see below), and $f$ is a source term. Here, boldface symbols are used for vector quantities. For a fixed time $t$, the PD Laplacian can be expressed as:

$$\mathcal{L}_\delta u(x) = \int_{H_x} \mu(\hat{x} - x)[u(\hat{x}) - u(x)]\mathrm{d}\hat{x} \tag{2}$$

where $\mu(x)$ is a non-negative even function, called *the kernel function* that defines the nonlocal interactions in neighborhood of spatial points [63-65]. In this work, we take $\mu$ to be an integrable function with compact support. Since $\mu(\hat{x} - x) = \mu(x - \hat{x})$, we have:

$$\mathcal{L}_\delta u(x) = \int_{H_x} \mu(x - \hat{x}) u(\hat{x})\mathrm{d}\hat{x} - u(x) \int_{H_x} \mu(x - \hat{x})\mathrm{d}\hat{x} \tag{3}$$

Assume $\mu(x)$ is defined over $(-\infty, +\infty)$, with $\mu=0$ outside of the horizon of $x=0$. With $\mu(x)$ being a given function, let $\beta = \int_{-\infty}^{+\infty} \mu(x)\mathrm{d}x$. The PD Laplacian becomes [64]:

$$\mathcal{L}_\delta u = \mu * u - \beta u \tag{4}$$

Where $(*)$ denotes the convolution integral operation.

### 2.1 Peridynamic Boundary Conditions

In problems specified by classical local theories, constraints are in the form of boundary conditions imposed on the surfaces of the 3D domain. In nonlocal problems, constraints are in the form of specified values on regions outside of the domain, where they have nonlocal interactions with parts of the domain [66]. Therefore, in the nonlocal problems *constrained-volume* and *volume-constraints* are used instead of boundaries and boundary conditions, respectively [66]. Such description of course, depends on the domain definition. For example, volume constraints may also be considered to be inside the domain. In this study however, the domain refers to the space where $u(x,t)$ is not specified and is solved for. Nevertheless, in many practical applications of peridynamics, imposing local-type boundary conditions is desired, for practical reasons.

Local boundary conditions can be enforced on a peridynamic body ($\Omega$), for example, via extending the domain by $\delta$ in the normal direction of the surface $\partial\Omega$. Quantities on the constrained volume, which is the domain extension $\Gamma$, are specified such that the local boundary condition is effectively reproduced on $\partial\Omega$ [67-70]. Values on $\Gamma$ are, in fact, volume constraints acting to enforce local boundary conditions. Fig. 2 schematically shows the peridynamic body $\Omega$, its boundary $\partial\Omega$, and the constrained-volume $\Gamma$.

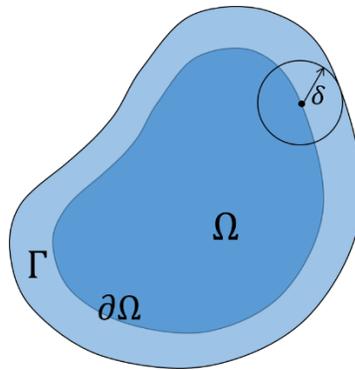

Fig. 2. Schematic of a peridynamic domain ($\Omega$), its boundary ($\partial\Omega$), and its constrained volume ($\Gamma$).

One way to impose volume constraints on $\Gamma$ with minimal or no difference from imposing local boundary conditions on $\partial\Omega$, is the scheme discussed in [67,68], known as the "fictitious nodes method". The terminology of this scheme refers to $\Gamma$ as a "fictitious" region since it is not a part of the domain. Note that in 2D and 3D this method will not be exact except for the simplest geometries [68]. In this scheme, the volume constraints are implicit and time-dependent, i.e. values on $\Gamma$ vary in time and are related to the values in the body $\Omega$ at that time. This type of implicit volume constraint can effectively impose a local BC on $\Omega$. The enforcement of some local BCs using the fictitious region method for the one-dimensional case follows.

To enforce the local Dirichlet BC:

$$u(\cdot,t)|_{\partial\Omega} = u(\delta,t) = u_b \tag{5}$$

with respect to the 1D configuration in Fig. 3, and a given $u_b$, values on $\Gamma$ should satisfy:

$$u(\cdot,t)|_{\Gamma} = u_\Gamma(x,t) = 2u_b - u(2\delta - x, t). \tag{6}$$

In order to apply the local Neumann BC:

$$\frac{\partial u}{\partial x}(\cdot, t)\Big|_{\partial \Omega} = \frac{\partial u}{\partial x}(\delta, t) = q_b, \quad (7)$$

given $q_b$, values in the constrained region are set as:

$$u(\cdot, t)|_\Gamma = u_\Gamma(x, t) = -2q_b(\delta - x) + u(2\delta - x, t) \quad (8)$$

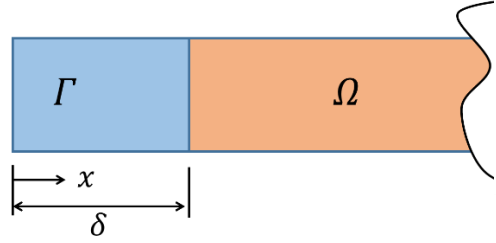

Fig. 3.  Schematic for the fictitious domain (constrained volume $\Gamma$) and the peridynamic body $\Omega$ in 1D. Time-dependent values on $\Gamma$ can be set to enforce some prescribed local boundary condition at $x = \delta$.

With this approach, similar to the Dirichlet BC, the Neumann BC is imposed by assigning value of $u$ rather than the values of its derivative. Note that Eq. (8) is new and different from the approach given in [67] where a source term is added in order to reproduce the Neumann BC. This specific form of Eq. (8) will show advantages when the spectral method will be used (see Section 3.2 below).

### 3. Spectral method for peridynamics with volume penalization
#### 3.1. Spectral method

Let $u(x, t)$ be a complex-valued function defined over the periodic domain $x \in T = [0, 2\pi]$ with 0 identified with $2\pi$, and evolve in time $t > 0$. Then $u(x, t)$ can be expressed with the infinite Fourier series in space:

$$u(x, t) = \sum_{k=-\infty}^{+\infty} \breve{u}_k(t) e^{\zeta k x}, \quad k = 0, \pm 1, \pm 2, \ldots \quad (9)$$

where $k$ is integer, $\zeta = \sqrt{-1}$, and:

$$\breve{u}_k(t) = \frac{1}{2\pi} \int_0^{2\pi} u(x, t) e^{-\zeta k x} \, dx, \quad k = 0, \pm 1, \pm 2, \ldots \quad (10)$$

are the Fourier coefficients of $u$ for different values of $k$. Eq. (10) is also called the Fourier transform of $u$ while Eq. (9) is the inverse Fourier transform relation.

Let the source term $f(x, t)$ and the kernel function with the form below, be also complex-valued functions defined over the periodic domain T.

$$\mu(x) = \begin{cases} \text{even function} & |x| \leq \delta \\ 0 & |x| > \delta \end{cases} \quad (11)$$

Then, the PD diffusion equation over the periodic domain T is:

$$\frac{\partial u}{\partial t} = \nu(\mu *_T u - \beta u) + f \tag{12}$$

where ($*_T$) denotes the "circular convolution" integral (aka "cyclic" or "periodic convolution") [71,72]:

$$\mu *_T u = \int_T \mu(x - \hat{x}) u(\hat{x}, t) \mathrm{d}\hat{x} \tag{13}$$

We approximate $u(x, t)$ by the *truncated (finite) Fourier series* of $u$:

$$u^N(x, t) = \sum_{k=-N/2}^{+N/2-1} \breve{u}_k(t) e^{\zeta k x} \tag{14}$$

Based on the finite Fourier series approximation, two similar, but not necessarily identical numerical schemes can be used for solving Eq. (12). One method is the *Fourier-Galerkin* method [73] in which the following weak form is solved:

$$\int_0^{2\pi} \left[ \frac{\partial u^N}{\partial t} - \nu(\mu^N *_T u^N) + \nu \beta u^N - f^N \right] e^{-\zeta k x} \mathrm{d}x = 0 \quad \text{for each } k = -\frac{N}{2}, \dots, \frac{N}{2} - 1 \tag{15}$$

Here $\mu^N$ and $f^N$ are the finite Fourier series approximations for $\mu$ and $f$. The integration on each term in Eq. (15) is the Fourier transform of that term. Eq. (15) is then equivalent to:

$$\frac{\mathrm{d}\breve{u}_k}{\mathrm{d}t} - 2\pi \nu \breve{\mu}_k \breve{u}_k + \nu \beta \breve{u}_k - \breve{f}_k = 0 \quad \text{for each } k = -\frac{N}{2}, \dots, \frac{N}{2} - 1 \tag{16}$$

We observe that, the circular convolution is transformed into a product operation of Fourier coefficients. In the Fourier-Galerkin method the ODE in Eq. (16) is solved for the Fourier coefficients of $u^N$. The solution can be transformed to the physical space with the inverse relation given by Eq. (14).

Another approach is the *Fourier Collocation* method [73] which is the one we will use in the present study. This method focuses on the solution in the physical space. The approximated solution in Eq. (14) is represented by its values at grid points $x_i = i\Delta x$, with $\Delta x = \frac{2\pi}{N}$ and $i \in \{0, \dots, N - 1\}$. In this method, $u^N(x, t)$ satisfies the strong form below at the collocation points $x_i$:

$$\frac{\partial u_i^N}{\partial t} - \nu\left(\mu_i^N *_T u_i^N\right) + \nu \beta u_i^N - f_i^N = 0 \tag{17}$$

Since $\mu$ and $u$ are approximated by finite Fourier series, the circular convolution can be evaluated by the inverse transform of the product of Fourier coefficients according to the convolution theorem [72]:

$$\mu^N *_T u^N = \mathcal{F}^{-1}(2\pi \breve{\mu}_k \breve{u}_k) \tag{18}$$

where $\mathcal{F}^{-1}$ refers to the inverse Fourier transform operation.

For practical applications of this method, a discrete-level operation is required to compute the Fourier transform and its inverse. The obvious choice is the Discrete Fourier Transform (DFT) [73]:

$$\tilde{u}_k(t) = \frac{1}{N} \sum_{i=0}^{N-1} u^N(x_i, t)\, e^{-\zeta k x_i} \tag{19}$$

and its inverse relation (iDFT)

$$u^N(x_i, t) = \sum_{k=-N/2}^{N/2-1} \tilde{u}_k(t)\, e^{\zeta k x_i} \tag{20}$$

Note that $\tilde{u}_k$ are approximations to the exact Fourier coefficients $\tilde{u}_k$.

Employing DFT, the Fourier Galerkin method yields to:

$$\frac{d\tilde{u}_k}{dt} - \nu(\widetilde{\mu^N *_T u^N}) + \nu\beta\tilde{u}_k - \tilde{f}_k = 0 \tag{21}$$

By the convolution theorem for DFT [71] we obtain:

$$\frac{d\tilde{u}_k}{dt} = \nu\tilde{\mu}_k \tilde{u}_k \Delta x - \nu\beta\tilde{u}_k + \tilde{f}_k \tag{22}$$

Using DFT and iDFT for transformation, the Fourier Collocation method in Eq. (17) becomes:

$$\frac{\partial u_i^N}{\partial t} = \nu \mathcal{F}_D^{-1}(\tilde{\mu}_k \tilde{u}_k \Delta x) - \nu\beta u_i^N + f_i^N \tag{23}$$

Where $\mathcal{F}_D^{-1}$ denotes the inverse DFT. Let $S$ denote the arc length of T (in our case $S = 2\pi$). If the periodic domain of computation is not $[0, S]$, i.e. it does not start at the origin at its left end, for example if it is $[-\frac{S}{2}, \frac{S}{2})$, then the kernel function may need to be shifted depending on the DFT solver (see appendix A).

The dominant computational cost in both methods is computing the DFT and its inverse, which are $O(N\log N)$ operations via FFT algorithms [45,46]. This is a significant improvement over the $O(N^2)$ cost for the meshfree collocation with one-point Gaussian quadrature or the FE methods used for PD problems. The extension of the spectral method to higher dimension is straightforward.

The above scheme works only for problems with periodic boundary conditions. We propose a penalization scheme that will allow us to apply spectral methods to general PD models with non-periodic boundary conditions in the next section.

### 3.2. Volume penalization

We employ the volume penalization (VP) technique developed for local problems in [55,56], to impose arbitrary volume constraints in a general PD problem.

In this method, the one-dimensional domain $\Omega$ is extended by $\delta$ at both ends as the constrained volume ($\Gamma$) to apply the nonlocal boundary conditions. The idea in the VP scheme is to consider periodicity for this extended domain, i.e. $T = \Omega \cup \Gamma$, and penalize the solution in the constrained domain to maintain the desired constraint values (See Fig. 4).

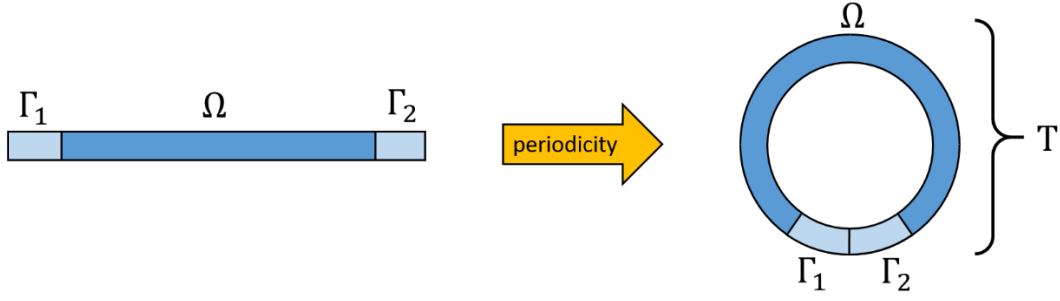

Fig. 4. Extension of 1D peridynamic non-periodic domain ($\Omega \cup \Gamma$) to the periodic domain T used in spectral method with volume penalization.

The PD diffusion equation is extended by adding a penalization term, which is zero on $\Omega$, but takes large values on $\Gamma$:

$$\frac{\partial u_\varepsilon(\boldsymbol{x},t)}{\partial t} = \nu \mathcal{L}_\delta u_\varepsilon(\boldsymbol{x},t) + f(\boldsymbol{x},t) - \frac{\chi(\boldsymbol{x},t)}{\varepsilon}[u_\varepsilon(\boldsymbol{x},t) - u_\Gamma(\boldsymbol{x},t)] \qquad (24)$$

In this equation $\varepsilon$ is a small number called here *the penalization factor*, $u_\varepsilon$ is the solution to the penalized PD diffusion equation, $u_\Gamma(\boldsymbol{x},t)$ is the volume constraint value at point $\boldsymbol{x} \in \Gamma$ and time $t$, and $\chi(\boldsymbol{x},t)$ is the following mask function:

$$\chi(\boldsymbol{x},t) = \begin{cases} 1 & \boldsymbol{x} \in \Gamma \\ 0 & \boldsymbol{x} \in \Omega \end{cases} \qquad (25)$$

For sufficiently small $\varepsilon$, the penalization term dominates on $\Gamma$:

$$\nu \mathcal{L}_\delta u_\varepsilon(\boldsymbol{x},t) + f(\boldsymbol{x},t) \ll \frac{1}{\varepsilon}[u_\varepsilon(\boldsymbol{x},t) - u_\Gamma(\boldsymbol{x},t)], \qquad (26)$$

leading to:

$$\frac{\partial u_\varepsilon(\boldsymbol{x},t)}{\partial t} \cong -\frac{1}{\varepsilon}[u_\varepsilon(\boldsymbol{x},t) - u_\Gamma(\boldsymbol{x},t)] \quad \text{on } \boldsymbol{x} \in \Gamma. \qquad (27)$$

Accordingly, this penalization term enforces an exponential decay for $u_\varepsilon$ to $u_\Gamma$ over the constrained domain. This effectively enforces the desired local boundary condition on $\partial\Omega$ if $u_\Gamma(\boldsymbol{x},t)$ is assigned via the scheme described in Section 2. Discussion on convergence of $u_\varepsilon$ to $u$ as $\varepsilon$ goes to zero is provided in Section 6.

To apply the spectral method, $u_\varepsilon$ is approximated with the finite Fourier series $u_\varepsilon^N$ on the periodic domain (T). To avoid complexity in notation, let $y = u_\varepsilon^N$. The spatially discretized version of Eq. (24) in 1D for using the *boundary-adapted spectral (BAS) method* is:

$$\frac{\partial y_i^n}{\partial t} = \nu \boldsymbol{\mathcal{F}}_D^{-1}(\tilde{u}_k \tilde{y}_k^n \Delta x) - \nu \beta y_i^n + f_i^n - \frac{\chi_i}{\varepsilon}(y_i^n - y_{\Gamma,i}^n) \qquad (28)$$

where the superscript $n$ refers to $n^{th}$ time step. A convergence study with respect to the spatial discretization size is provided in Section 6.

Any applicable temporal integration scheme maybe used to now solve the first order ODE in Eq. (28), and update the solution at each time step. With the Forward Euler method, for example, we have:

$$y_i^{n+1} \approx y_i^n + \Delta t \frac{dy_i^n}{dt} \tag{29}$$

where $\Delta t$ is the time step. The stability restriction on the time step size for this explicit method is derived in the next section.

### 4. Stability analysis

Here we follow the stability analysis in [35] to find the restriction on time steps for the BAS method with VP using the explicit Euler time integration scheme. It can be shown that Eq. (28) in the physical space is algebraically equivalent to:

$$\frac{dy_i^n}{dt} = \nu \sum_{j=0}^{N-1} \mu_{i-j} y_j^n \Delta x - \nu \beta y_i^n + f_i^n - \frac{\chi_i}{\varepsilon}(y_i^n - y_{\Gamma,i}^n) \tag{30}$$

where $\mu_{i-j} = \mu(x_i - x_j)$. Note that in the term in the summation above is zero for all $x_j \notin [x_i - \delta, x_i + \delta]$. Although similar, the discretized volume integral in Eq. (30) is not identical to one-point Gaussian quadrature used in the meshfree method. Functions are approximated with the truncated Fourier series, which is not the case in the conventional meshfree method.

With Forward Euler (first-order explicit) temporal integration, Eq. (30) becomes:

$$\frac{y_i^{n+1} - y_i^n}{\Delta t} = \nu \sum_j \mu_{i-j} y_j^n \Delta x - \nu \beta y_i^n + f_i^n - \frac{\chi_i}{\varepsilon}(y_i^n - y_{\Gamma,i}^n) \tag{31}$$

Take:

$$y_i^n = \lambda^n e^{\zeta k x_i} \tag{32}$$

where $\lambda$ is a complex number. Substituting Eq. (32) to Eq. (31), results in:

$$\frac{\lambda^{n+1} - \lambda^n}{\Delta t} e^{\zeta k x_i} = \nu \sum_j \mu_{i-j} \lambda^n e^{\zeta k x_j} \Delta x - \nu \beta \lambda^n e^{\zeta k x_i} + f_i^n - \frac{\chi_i}{\varepsilon}(\lambda^n e^{\zeta k x_i} - u_{\Gamma,i}^n) \tag{33}$$

For simplicity, let $u_{\Gamma,i}^n = 0$, $f_i^n = 0$, and $\rho = \frac{\lambda^{n+1}}{\lambda^n}$ for every $i$ and $n$, then:

$$\frac{\rho - 1}{\Delta t} = \nu \sum_j \mu_{i-j} e^{\zeta k(x_j - x_i)} \Delta x - \nu \beta - \frac{\chi_i}{\varepsilon} \tag{34}$$

Since $\mu_{i-j} = \mu_{j-i}$, let $p = j - i$, and $x_p = x_j - x_i$, we obtain:

$$\frac{\rho - 1}{\Delta t} = \nu \sum_p \mu_p e^{\zeta k x_p} \Delta x - \nu \beta - \frac{\chi_i}{\varepsilon} \tag{35}$$

If $m = \frac{\delta}{\Delta x}$ is integer, then, since $\mu_p = \mu_{-p}$ and $x_{-p} = -x_p$:

$$\sum_{p=-m}^{+m} \mu_p e^{\zeta k x_p} \Delta x = \sum_{p=1}^{m}\left(\mu_p e^{\zeta k x_p} \Delta x + \mu_{-p} e^{\zeta k x_{-p}} \Delta x\right) + \mu_0 \Delta x = 2 \sum_{p=0}^{m} \mu_p \cos(k x_p) \Delta x \tag{36}$$

$$= \sum_{p=-m}^{m} \mu_p \cos(k x_p) \Delta x$$

Substituting Eq. (36) into Eq. (35) results in:

$$\frac{\rho - 1}{\Delta t} = \nu \left( \sum_{p=-m}^{m} \mu_p \cos(kx_p) \Delta x - \beta \right) - \frac{\chi_i}{\varepsilon} \tag{37}$$

Define $M$:

$$M = \nu \left( \sum_{p=-m}^{m} \mu_p \cos(kx_p) \Delta x - \beta \right) \tag{38}$$

Then:

$$\rho = \left( M - \frac{\chi_i}{\varepsilon} \right) \Delta t + 1 \tag{39}$$

To maintain stability, we seek $\Delta t$ such that $|\rho| \leq 1$. Therefore:

$$\left| \left( M - \frac{\chi_i}{\varepsilon} \right) \Delta t + 1 \right| \leq 1 \tag{40}$$

or,

$$-1 \leq \left( M - \frac{\chi_i}{\varepsilon} \right) \Delta t + 1 \leq 1 \tag{41}$$

or equivalently:

$$0 \leq \left( \frac{\chi_i}{\varepsilon} - M \right) \Delta t \leq 2 \tag{42}$$

In order to satisfy the left inequality in Eq. (42), since $\chi_i$ is either 0 or 1, and that $\varepsilon$ and $\Delta t$ are positive quantities, we need $M \leq 0$. According to Eq. (38), this imposes the following condition on $\mu(x)$:

$$\sum_{p=-m}^{m} \mu_p \cos(kx_p) \Delta x \leq \beta = \int_{-\delta}^{\delta} \mu(x) \mathrm{d}x \tag{43}$$

This condition holds for sufficiently small $\Delta x$, since:

$$\int_{-\delta}^{\delta} \mu(x) \cos(kx) \, \mathrm{d}x \leq \int_{-\delta}^{\delta} \mu(x) \mathrm{d}x \tag{44}$$

Most kernel functions in use satisfy the requirement in Eq. (43).

The inequality in Eq. (42) also requires:

$$\Delta t \leq \frac{2}{\frac{\chi_i}{\varepsilon} - M} \tag{45}$$

From Eqs. (38) and (43):

$$-M \leq 2\nu\beta \tag{46}$$

According to Eqs. (45) and (46), and that $\chi_i$ is either 0 or 1, the following restriction on $\Delta t$ for stable solution is suggested:

$$\Delta t \leq \frac{2}{\frac{1}{\varepsilon}+2\nu\beta}. \tag{47}$$

Even if the restriction above is obtained assuming zero values for $u_\Gamma$ and $f$, we find that it is also sufficient for obtaining stable results for the examples shown in Section 5, with nonzero $u_\Gamma$ and $f$.

From Eq. (47), we see that penalization puts a stronger restriction on the time step compared with the condition found with the conventional meshfree method [67]. However, the cost due to the increased number of time steps is likely to be overcome by the gains in the complexity order (in terms of node number) when computing the convolution integral with FFT. Note that the increased number time steps does affect the complexity order, since according to Eq. (47), the time step size does not depend on $N$.

### 5. Example problems and discussion

We now compare the performance of the PD spectral method with the regular integration (one point Gaussian quadrature) of the convolution integral in the PD Laplacian. Then we analyze two one-dimensional nonlocal diffusion problems to demonstrate the capability of the BAS method introduced. The first problem has local Dirichlet boundary conditions at both ends, while the second has local Neumann boundary conditions at both ends.

For these examples we need to select a kernel function $\mu(x)$. According to [74], one possible choice for the kernel function is of the form:

$$\mu(\boldsymbol{x}) = \frac{(4-\alpha)(3-\alpha)}{\delta^{(3-\alpha)}}\left(1 - \frac{|x|}{\delta}\right)\frac{1}{|x|^\alpha} \tag{48}$$

where $\alpha$ can take values 0, 1, and 2, for example (for details see [74]). For the case $\alpha = 0$, we have:

$$\mu(x) = \frac{12}{\delta^3}\left(1 - \frac{|x|}{\delta}\right) \tag{49}$$

Chen and Bobaru [74] showed that a constructive approach to a peridynamic kernel leads to the choice of $\alpha = 2$. Here, however, we choose $\alpha = 0$ for simplicity and to avoid the singularity when calculating $\beta = \int_{H_0} \mu(x) dx$. For the case $\alpha = 0$, this integral is $\frac{12}{\delta^2}$. For the other values of $\alpha$, $\beta$ can be calculated using the Cauchy principal value [75].

### 5.1. Efficiency of the peridynamic spectral method

Here we compute $\mathcal{L}_\delta u$ for $u = \sin(\pi x)$ with $\delta = 0.2$ in $x \in [-1,1]$ via two methods: the direct numerical integration using one-point Gaussian quadrature:

$$(\mathcal{L}_\delta u)_i = \sum_{j=i-\text{round}(\delta/\Delta x)}^{i+\text{round}(\delta/\Delta x)} \mu_{i-j} u_j \Delta x - \beta \tag{50}$$

and the spectral method:

$$(\mathcal{L}_\delta u)_i = \mathcal{F}_D^{-1}(\tilde{\mu}_k \tilde{u}_k \Delta x) - \beta \tag{51}$$

The kernel function $\mu$ is the one defined by Eq.(49), and therefore $\beta = \frac{12}{\delta^2}$.

$\mathcal{L}_\delta u$ is computed using both methods for several discretization sizes, with $N$ varying between $2^8$ to $2^{20}$. Computation are performed using *MATLAB 2018a* on a Dell-Precision T7810 workstation PC, with twenty logical Intel(R) Xeon(R) CPU E5-2687 W v3 @3.10 GHz processors, and 64 GB of installed memory. To provide a fair comparison between methods, we restricted MATLAB to use only one computational thread, to prevent multi-core processing with its FFT solver. The computational time for each method to calculate $\mathcal{L}_\delta u$ for various discretization sizes are provided in Table 1.

Table 1. Comparison of run-times between the one-point Gaussian quadrature and the spectral method in calculating the peridynamic Laplacian

| $N$ (number of nodes) | Gaussian quadrature time (sec) | Spectral method time (sec) |
|---|---|---|
| $2^8$ = 256 | 1.79e-4 | 2.11e-5 |
| $2^{12}$ = 4,096 | 5.03e-2 | 1.63e-4 |
| $2^{16}$ = 65,536 | 1.07e+1 | 5.04e-3 |
| $2^{20}$ = 1,048,576 | 3.00e+3 (50 min) | 6.73e-2 (67 ms) |

As observed, the time for one-point Gaussian quadrature is $O(N^2)$, while the spectral method performs even more efficient than $O(N\log N)$. The reason for the over-performance of the spectral method may be due to the efficient FFT solver in MATLAB, which uses optimized algorithms with respect to data size and structure.

### 5.2. Diffusion with Dirichlet boundary condition

We now solve an example of PD diffusion problem with local Dirichlet BCs, using the BAS method, and compare the numerical solution with the analytical solution. Consider, for example:

$$u(x,t) = \frac{2x}{L} + e^{-vt} \sin\left(\frac{2\pi x}{L}\right) \tag{52}$$

The function in Eq. (52) is the exact solution to the following nonlocal diffusion problem over the domain $\Omega = \left[-\frac{L}{2}, \frac{L}{2}\right]$:

$$\frac{\partial u(x,t)}{\partial t} = v \mathcal{L}_\delta u(x,t) + f(x,t), \tag{53}$$

with

$$f(x,t) = v\left\{\frac{6L^2}{\delta^4 \pi^2}\left[\cos\left(\frac{2\pi\delta}{L}\right) - 1\right] + \frac{12}{\delta^2} - 1\right\} e^{-vt} \sin\left(\frac{2\pi x}{L}\right), \tag{54}$$

the initial condition:

$$u(x,0) = \frac{2x}{L} + \sin\left(\frac{2\pi x}{L}\right), \tag{55}$$

and the local Dirichlet boundary conditions:

$$u\left(-\frac{L}{2}, t\right) = -1 \tag{56}$$

$$u\left(\frac{L}{2}, t\right) = 1 \tag{57}$$

The manufactured solution $u(x,t)$ in Eq. (52) has a special property: with the fictitious nodes scheme described in Section 2, the values of $u$ in $\left[-\frac{L}{2} - \delta, -\frac{L}{2}\right)$ and $\left(\frac{L}{2}, \frac{L}{2} + \delta\right]$, satisfy the volume constraint relationship in Eq. (6). This property then makes is easier to find the $f(x,t)$ above. Without this property, one can still find $f(x,t)$ analytically, but as a relatively more complex piecewise function.

To solve this problem with the proposed method, we select $L = 2$, $v = 0.2$, $\delta = 0.2$, and the total diffusion time $t_{\max} = 15$. The computational domain is then extended to $T = \Omega \cup \Gamma : \left[-\frac{L}{2} - \delta, \frac{L}{2} + \delta\right) = [-2.2, 2.2)$, with $[-2.2, -2)$ and $(2, 2.2)$ being the constrained domains $\Gamma_1$ and $\Gamma_2$ respectively. Note that the computational domain interval does not include $x = 2.2$ due to the periodicity of the spectral method that mandates $x = 2.2$ be identical to $x = -2.2$. Choosing the number of spatial nodes to be a power of two ($N = 2^P$) has certain benefits in parallelization of FFT algorithms [46]. The extended domain $[-2.2, 2.2)$ is discretized with $N = 2^9$ nodes. The time step and penalization factor are selected as $\Delta t = 5 \times 10^{-4}$ and $\varepsilon = 5 \times 10^{-4}$, respectively. The algorithm for the implementation of the proposed method is provided in the appendix A. According to the fictitious nodes scheme described in Section 2, volume constraint values on $\Gamma$ are calculated explicitly from Eqs. (58) and (59), using the solution on $\Omega$ at the previous time step. For $x \in \Gamma_1 : \left[-\frac{L}{2} - \delta, -\frac{L}{2}\right)$:

$$y_{\Gamma 1}(x_i, t^{n+1}) = 2u_{b1} - y(2\delta - x_i, t^n), \tag{58}$$

where $u_{b1}$ is the Dirichlet BC value given in Eq. (56), and $y$ is $u_\varepsilon^N$: the numerical solution using the BAS method with VP. For $x \in \Gamma_2 : \left(\frac{L}{2}, \frac{L}{2} + \delta\right)$, the equation below applies the boundary condition:

$$y_{\Gamma 2}(x_i, t^{n+1}) = 2u_{b2} - y(2L + 2\delta - x_i, t^n), \tag{59}$$

where $u_{b2}$ is the BC value given in eq. (57). Note that the volume constraints can be applied implicitly as well by replacing $t^n$ with $t^{n+1}$ in Eq. (58) and (59), where one need to iterate until $y_\Gamma^{n+1}$ converges.

Fig. 5 shows the time snapshots of the nonlocal diffusion process. We observe the excellent match with the analytical nonlocal solution for this non-periodic problem.

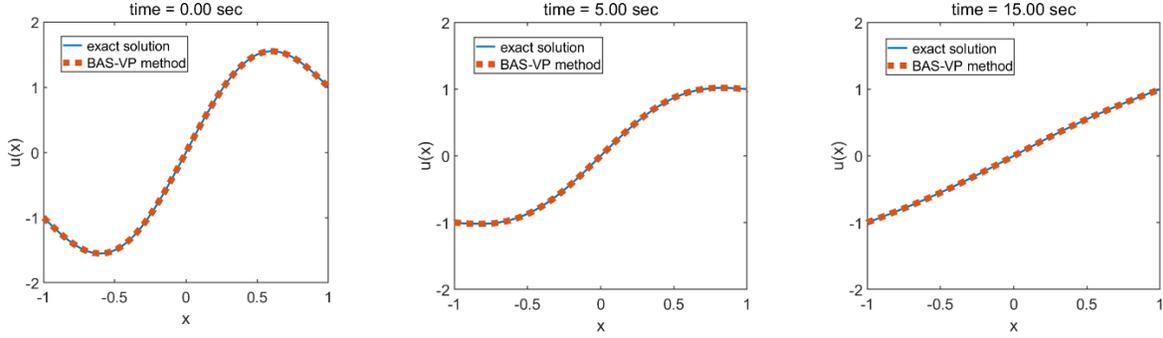

Fig. 5.  Comparison between the analytical solution of a nonlocal 1D diffusion process, with non-periodic (Dirichlet) boundary conditions, and the solution obtained by the peridynamic spectral method with volume penalization.

The absolute error distribution, normalized by the infinity norm of the initial data function ($\frac{|u-u_\varepsilon^N|}{\|u_0\|_\infty}$) is plotted in Fig. 6. We observe the rise and decay of the error in time, in the interior region of the domain $\Omega$, while the error near the boundaries rises and approaches to permanent amount (see Video 1).

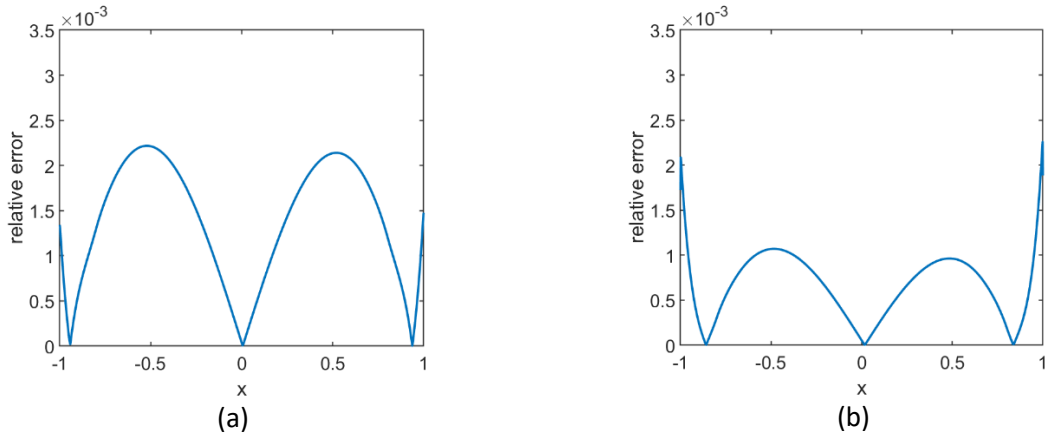

(a)                                                                 (b)

Fig. 6.  Time snapshots of the relative error in 1D nonlocal diffusion with non-periodic (Dirichlet) boundary conditions using the peridynamic spectral method with volume penalization. a) $t = 5$; b) $t = 15$.

The maximum relative error ($\frac{\|u-u_\varepsilon^N\|_\infty}{\|u_0\|_\infty}$) is plotted versus diffusion time in Fig. 7. The slope discontinuity in this plot can be understood by observing the behavior in Video 1: the location of the maximum relative error switches to the boundaries, as time progresses.

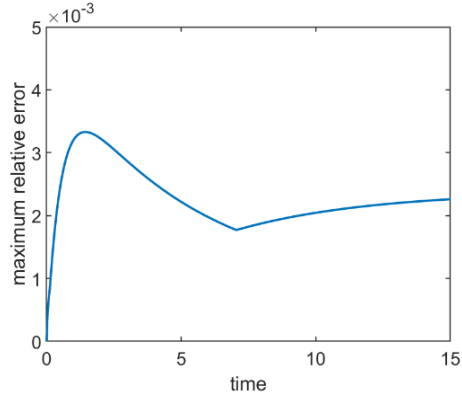

Fig. 7. Variation of the maximum relative error in time for the 1D nonlocal diffusion example with non-periodic (Dirichlet) boundary conditions, using the peridynamic spectral method with volume penalization.

This plot suggests that the observed error time-evolution consists of two periods: at first, the error in the domain interior dominates, while later the error near the boundaries becomes more important. Our parametric studies in appendix B show that the decaying error on the interior originates from the spatial discretization and is reduced with grid refinement. The error near the boundaries, however, depends of the penalization factor and is reduced by selecting a smaller $\varepsilon$. Convergence studies for the total error in terms of discretization size and penalization parameter are given in Section 6.

### 5.3. Diffusion with Neumann boundary condition

To demonstrate the capability of the proposed BAS method in solving PD problems with arbitrary boundary conditions we now discuss and example with Neumann BCs. Consider the function:

$$u(x,t) = \frac{2x}{L} + e^{-vt} \cos\left(\frac{2\pi x}{L}\right), \tag{60}$$

the exact solution for the nonlocal diffusion equation in Eq. (53) on the interval $\left[-\frac{L}{2}, \frac{L}{2}\right]$, with

$$f(x,t) = v\left\{\frac{6L^2}{\delta^4\pi^2}\left[\cos\left(\frac{2\pi\delta}{L}\right) - 1\right] + \frac{12}{\delta^2} - 1\right\} e^{-vt} \cos\left(\frac{2\pi x}{L}\right), \tag{61}$$

initial condition

$$u(x,0) = \frac{2x}{L} + \cos\left(\frac{2\pi x}{L}\right), \tag{62}$$

and local Neumann BCs:

$$\frac{\partial u}{\partial x}\left(-\frac{L}{2}, t\right) = 1 \tag{63}$$

$$\frac{\partial u}{\partial x}\left(\frac{L}{2}, t\right) = 1. \tag{64}$$

Similar to the previous example for Dirichlet BC, the manufactured solution in Eq. (60) satisfies the volume constraint relationship in Eq. (8) when $x \in \Gamma_1$ or $\Gamma_2$.

For this problem, $L$, $\delta, \nu$, and $t_{max}$, are the same as in the previous example in Section 5.1. We use the explicit implementation of the fictitious nodes scheme for applying the volume constraints corresponding to Neumann BC (see Eq.(8)). For $x \in \Gamma_1: \left[-\frac{L}{2}-\delta, -\frac{L}{2}\right)$ we get:

$$y_{\Gamma_1}(x_i, t^{n+1}) = -2q_{b1}(\delta - x_i) + y(2\delta - x_i, t^n), \tag{65}$$

while for $x \in \Gamma_2: \left(\frac{L}{2}, \frac{L}{2}+\delta\right)$, we have:

$$y_{\Gamma_2}(x_i, t^{n+1}) = -2q_{b2} + y(2L + 2\delta - x_i, t^n) \tag{66}$$

Values for both $q_{b1}$ and $q_{b2}$ are 1 according to Eqs. (63) and (64). For the numerical solution, $N, \Delta t$, and $\varepsilon$ are the same as in the previous example.

Fig. 8 shows the evolution in time of the numerical solution $y$ in comparison with the exact nonlocal solution $u$. The results support the fact that the BAS method is capable solve peridynamic problems with arbitrary boundary conditions.

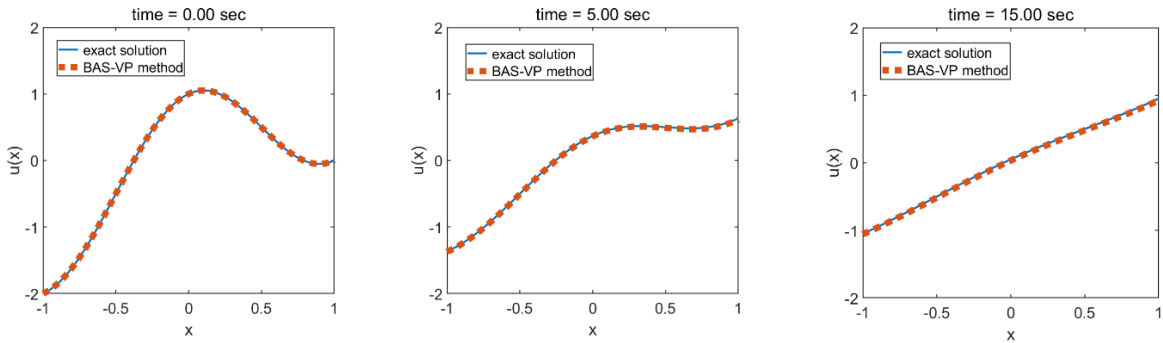

Fig. 8. Comparison between the analytical solution of a nonlocal 1D diffusion process, with non-periodic Neumann boundary conditions, and the solution obtained by the peridynamic spectral method with volume penalization.

## 6. Convergence

In this section we first provide a brief background on error estimates for the volume penalization method, and then present convergence studies with respect to the penalization and spatial discretization for the example problem shown in Section 5.2.

Angot et al.[55] proved that the solution of Navier-Stokes equation with the volume penalization in a periodic domain converges to the solution of Navier-Stokes equation with the proper exact boundary conditions, as $\varepsilon$ goes to zero. The error in the main domain for that case is shown to be at most of the $O(\varepsilon^{3/4})$. In the case of the classical diffusion equation, Kevlahan and Ghidaglia [56] showed for a specific problem that the error between the penalized periodic solution and the exact solution to the diffusion equation with non-periodic boundary conditions is at most $O(\varepsilon^{1/2})$. They observed that the computed error is $O(\varepsilon)$.

While a rigorous mathematical convergence analysis and error estimate for nonlocal diffusion equation with penalization and spectral method would be ideal, in this study we only provide some numerical results. The complete theoretical analysis is left for the future. We perform convergence studies on the example with two Dirichlet BCs discussed in Section 5.2.

First we study the convergence of the penalized periodic solution $u_\varepsilon^N$ (the solution to Eq. (24)), to the exact solution (Eq. (52)) of the un-penalized diffusion equation with Dirichlet BC. To this aim, we need to choose a relatively large $N$, and relatively small time step and keep them fixed while decreasing $\varepsilon$ in each test. This makes that the discretization and temporal integration errors minor compared with the penalization error which we want to investigate. We need to also choose a time span to approach the steady state where according to the observations in Fig. 6 and Fig. 7 permanent penalization error is dominated and remains relatively constant. The selected parameters for this convergence study are $N = 2^{15}$, $\Delta t = 1.97 \times 10^{-4}$, and $t_{max} = 30$, while $\varepsilon^{-1}$ value varies for each test. $\Delta t$ satisfies Eq. (47) restriction with the smallest penalization factor used in this convergence study ($\varepsilon = 1 \times 10^{-4}$). Fig. 9 shows the relative error versus $\varepsilon^{-1}$ for each test.

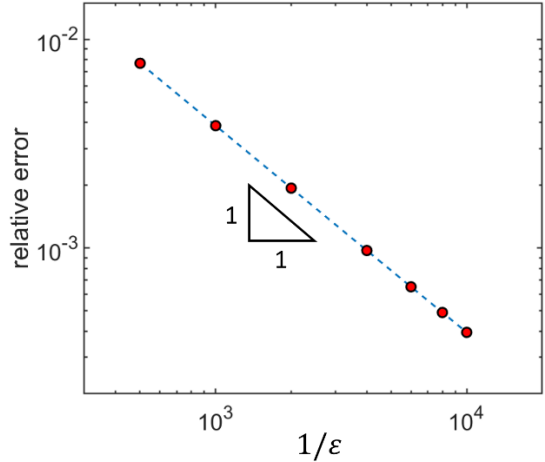

Fig. 9. Convergence study in terms of penalization factor for the peridynamic spectral method with volume penalization on the problem with Dirichlet boundary conditions shown in Section 5.2.

The results show that the penalization error varies with $O(\varepsilon)$, which is consistent with observations for the penalized classical diffusion equation [56].

To observe the convergence behavior with respect to spatial discretization size, we compare the maximum error in the whole time span for various $N$ values, while keeping constant the relatively small values of $\Delta t$ and $\varepsilon$. To this aim, we obtained the error $\max\limits_{0<t<t_{max}} \frac{\|u - u_\varepsilon^N\|_\infty}{\|u_0\|_\infty}$, with respect to different $N$ values in five tests with $\varepsilon = 5 \times 10^{-6}$, and $\Delta t = 9.99 \times 10^{-6}$. Again, $\Delta t$ satisfies the stability condition given in Eq. (47). Results are plotted in Fig. 10.

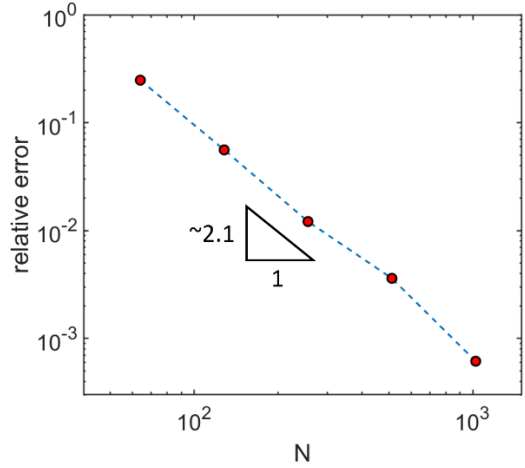

Fig. 10.  Convergence study with respect to the spatial discretization size for the spectral method with volume penalization on the problem with Dirichlet boundary conditions shown in Section 5.2.

As observed, the spatial convergence rate of peridynamic BAS method is $O(\Delta x^2)$ for this example problem. This is similar to meshfree-collocation method with one-point Gaussian quadrature [35].

Note that the all the comparison of the solutions with exact solution for obtaining errors are considered within the domain of interest $\Omega = T \backslash \Gamma$ which disregards the solution values on the penalized region $\Gamma$.

The general error of the peridynamic BAS method is bounded by the summation of the penalization error, spatial discretization errors (finite Fourier series approximation and DFT), and the explicit time integration error which for Forward Euler is known to be $O(\Delta t)$.

For the presented example with two Dirichlet BCs, the error of the introduced method appears to be bounded by $O(\varepsilon) + O(\Delta x^2) + O(\Delta t)$.

## 7. Conclusions

In this study we introduced a boundary-adapted spectral (BAS) method for peridynamic (PD) diffusion problems with arbitrary boundary conditions. The spectral approach transforms the convolution integral into a multiplication in the Fourier space, resulting in computations that scale as $O(N\log N)$. We demonstrated the efficiency of this method by comparing it with the commonly used one-point Gaussian quadrature method for spatial integration in a peridynamic model. In 1D, a problem with roughly one million nodes is solved in milliseconds with the spectral method whereas the one-point Gaussian quadrature approach required close to one hour in our tests.

A stability analysis for the peridynamic BAS method (with the volume-penalization approach) with Forward-Euler time integration for peridynamic transient diffusion problems suggested that the restriction on the time-step varies linearly with the penalization factor, for a sufficiently large one. We examined the performance of the method introduced for arbitrary boundary conditions with two examples of peridynamic transient diffusion using local Dirichlet and Neumann boundary conditions. We compared our numerical results against exact nonlocal solutions. Our convergence studies show that the error scales linearly with the penalization factor and quadratically with the discretization size. The method can be easily extended to other peridynamic/nonlocal models, in 1D and in higher dimensions.

**Acknowledgments**

This work has been supported by the AFOSR MURI Center for Materials Failure Prediction through Peridynamics (program managers Jaimie Tiley, David Stargel, Ali Sayir, Fariba Fahroo, and James Fillerup), by the ONR project #N00014-15-1-2034 "SCC: the Importance of Damage Evolution in the Layer Affected by Corrosion" (program manager William Nickerson), and by a Nebraska System Science award. The research of A.L. was supported in part by the NSF grant no. DMS-1716801.

**Appendix A: Boundary-adapted spectral method implementation for PD diffusion in MATLAB**

Here the MATLAB implantation for boundary-adapted spectral method with volume penalization (BAS-VP) is provided. First, note that Eq. (23) can be directly use when the periodic domain of computation is $[0, S)$, meaning the origin locates on the left end of the domain. If the domain of choice is $[b, b + S)$ then the following modified form of Eq. (23) should be used:

$$\frac{du_i^N}{dt} = v\mathcal{F}_D^{-1}\left(\widetilde{\mu^s}_k \tilde{u}_k \Delta x\right) - v\beta u_i^N + f_i^N \tag{67}$$

where $\widetilde{\mu^s}_k$ is the DFT of the shifted kernel function:

$$\mu^s(x) = \mu(x - b) \tag{68}$$

The reason is that the DFT definitions that governs the FFT solvers are based on $[0, S)$ domain. If $b = 0$, then the kernel function does not shift and Eq. (67) becomes identical to Eq. (23).

A MATLAB implementation of the peridynamic BAS method with VP for the transient diffusion example in Section 5.2 is as follows:

- Inputs:
    - Physical parameters: $v$, $\delta$, $f(x,t)$, $\mu(x)$, $L$, $t_{\max}$
    - Initial and boundary conditions: $u(x, 0)$, $u\left(-\frac{L}{2}, 0\right) = u_{b1}$, $u\left(\frac{L}{2}, 0\right) = u_{b2}$
    - BAS method with VP parameters: $N$, $\varepsilon$
- Initialization:
    - Calculate grid size: $\Delta x = \frac{L+2\delta}{N}$ (length of the extended domain divided by $N$)
    - Calculate time step: $\Delta t$ from Eq. (47) with $v$, $\varepsilon$, and $\mu(x)$
    - Discretize the extended domain: $x_i = -\frac{L}{2} - \delta + (i-1)\Delta x$ and $i = 1,2,...,N$
    - Shift the kernel function based on left-end of the extended domain and discretize: $\mu_i^s = \mu\left(x_i + \frac{L}{2} + \delta\right)$
    - Discretize the initial condition and the Source term: $y_i^0 = u(x_i, 0)$; $f_i^0 = f(x_i, 0)$;
    - Fast Fourier transform $\mu_i^s$ and $y_i^0$: $\widetilde{\mu^s}_k = \mathbf{FFT}(\mu_i^s)$ and $\widetilde{y^0}_k = \mathbf{FFT}(y_i^0)$

- Define constrained regions: $\Gamma_1 = x_i \in \left[-\frac{L}{2} - \delta, -\frac{L}{2}\right)$ and $\Gamma_2 = x_i \in \left(\frac{L}{2}, \frac{L}{2} + \delta\right)$
- Define the main domain: $\Omega = x_i \in \left[-\frac{L}{2}, \frac{L}{2}\right]$
- Discretize the mask function: $\chi_i = \chi(x_i)$ from Eq. (25).
- Calculate volume constraints on $\Gamma_1$ and $\Gamma_2$ from Eq. (58) and (59): $u_{\Gamma 1}(\Gamma_1, 0), u_{\Gamma 2}(\Gamma_2, 0)$
- Define $y_{\Gamma i}^0 = \begin{cases} \text{Eq. (58)} & x_i \in \left[-\frac{L}{2} - \delta, -\frac{L}{2}\right) \\ 0 & x_i \in \left[-\frac{L}{2}, \frac{L}{2}\right] \\ \text{Eq. (59)} & x_i \in \left(\frac{L}{2}, \frac{L}{2} + \delta\right) \end{cases}$
- Initialize step counter: $n = 0$
- Initialize time: $t^n = 0$

■ Solve the transient diffusion: while $t^n < t_{\max}$
- Update time: $t^{n+1} = t^n + \Delta t$
- Update solution: $y_i^{n+1} = y_i^n + \Delta t \left[ \nu \mathbf{FFT}^{-1}(\widetilde{\mu^s}_k \widetilde{y^n}_k \Delta x) - \nu \beta y_i^n + f_i^n - \frac{\chi_i}{\varepsilon}(y_i^n - y_{\Gamma,i}^n) \right]$
- Update the source term: $f_i^{n+1} = f_i(x_i, t^{n+1})$
- Update volume constraints: $y_{\Gamma i}^{n+1} = \begin{cases} \text{Eq. (58)} & x_i \in \left[-\frac{L}{2} - \delta, -\frac{L}{2}\right) \\ 0 & x_i \in \left[-\frac{L}{2}, \frac{L}{2}\right] \\ \text{Eq. (59)} & x_i \in \left(\frac{L}{2}, \frac{L}{2} + \delta\right) \end{cases}$
- Fast Fourier transform $y_i^{n+1}$: $\widetilde{y^{n+1}}_k = \mathbf{FFT}(y_i^{n+1})$
- Update step counter: $n = n + 1$

The algorithm above is for the example with Dirichlet BCs. In the case of Neumann BCs Eq. (58) and (59) are replaced with Eq. (65) and (66).

**Appendix B: Discretization error versus penalization error in BAS-VP method**

To obtain a better understanding of error distribution on the domain for the example in Section 5.2, and the evolution of maximum error during the diffusion process (see Fig. 6 and Fig. 7), we conducted two more simulations: One simulation with a much smaller $\varepsilon$ compared with the test Section 5.2, but the same $N$; and one simulation with a much larger $N$ compared to that test, but the same $\varepsilon$. The first simulation reveals the error behavior with respect to the discretization, while the second one is focused on the penalization error.

Results for the first simulation with $\varepsilon = 5 \times 10^{-6}$, $N = 2^9$ are given in Fig A.1.

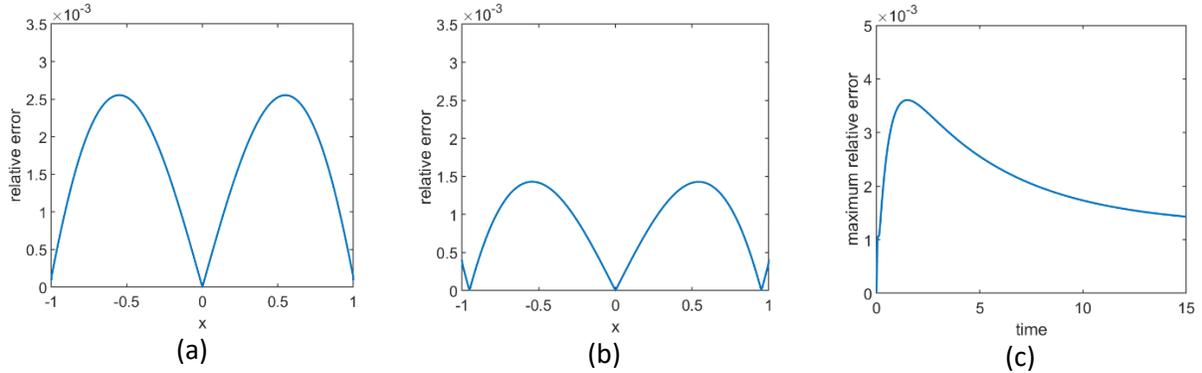

Fig A.1.  Time snapshots of the relative error dominated by discretization in 1D nonlocal diffusion problem in Section 5.2. at $t = 5$ (a), and $t = 15$ (b). The time-evolution of the maximum error in (c).

The second simulation is performed with $\varepsilon = 5 \times 10^{-4}$ and $N = 2^{15}$. Results are given in Fig A.2.

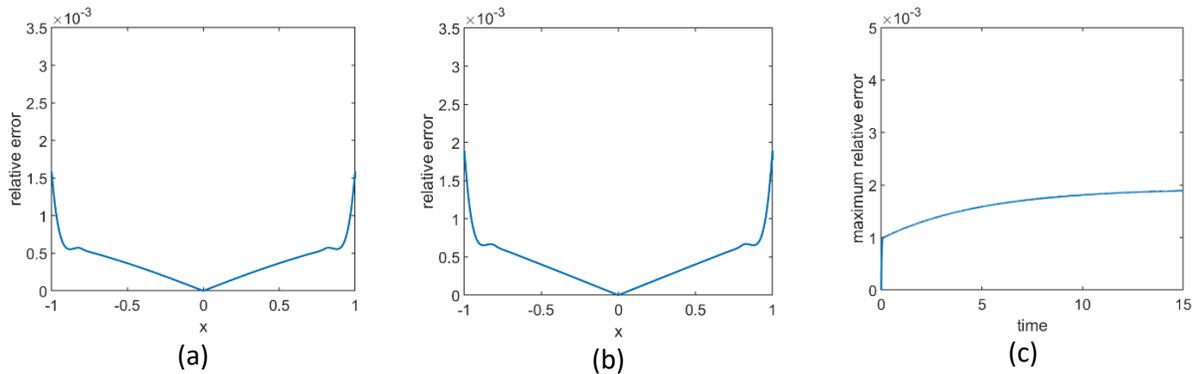

Fig A.2.  Time snapshots of the relative error dominated by penalization in 1D nonlocal diffusion problem shown in Section 5.2. at $t = 5$ (a), and $t = 15$ (b). The time-evolution of maximum error in (c).

The discretization error rapidly grows and then decays, while the penalization error grows near the boundaries and approaches a constant value in time. Comparing Fig A.1 and Fig A.2, with Fig. 6 and Fig. 7 in Section 5.2 (see also video 1), helps us to clearly identify the "mixture" of the penalization and discretization errors in the example corresponding to Fig. 6 and Fig. 7.